\documentclass[11pt,reqno]{amsart}

\usepackage[T1]{fontenc}
\usepackage[utf8]{inputenc}
\usepackage{lmodern}
\usepackage[margin=1in]{geometry}
\usepackage{amsmath,amssymb,mathtools}
\usepackage{enumitem}
\usepackage{microtype}
\usepackage{needspace}
\usepackage{cite}
\usepackage[hidelinks]{hyperref}

\newtheorem{theorem}{Theorem}[section]
\newtheorem{proposition}[theorem]{Proposition}
\newtheorem{lemma}[theorem]{Lemma}
\newtheorem{corollary}[theorem]{Corollary}
\newtheorem{question}[theorem]{Question}
\newtheorem{problem}[theorem]{Problem}
\theoremstyle{definition}
\newtheorem{definition}[theorem]{Definition}
\theoremstyle{remark}
\newtheorem{remark}[theorem]{Remark}

\DeclareMathOperator{\CA}{CA}
\DeclareMathOperator{\ICA}{ICA}
\DeclareMathOperator{\End}{End}
\DeclareMathOperator{\Fix}{Fix}
\DeclareMathOperator{\supp}{supp}
\DeclareMathOperator{\im}{Im}
\newcommand{\cL}{\mathcal{L}}
\newcommand{\id}{\mathrm{id}}
\newcommand{\one}{\mathbf{1}}

\title{When do lazy and invertible cellular automata generate all cellular automata?}

\author{Edgar Alcalá-Arroyo}
\address{Universidad de Guadalajara, Centro Universitario de Ciencias Exactas e Ingenierías, Departamento de Matem\'aticas. Guadalajara, México.}
\email{edgar.alcala7434@alumnos.udg.mx}

\author{Alonso Castillo-Ramirez}
\address{Universidad de Guadalajara, Centro Universitario de Ciencias Exactas e Ingenierías, Departamento de Matem\'aticas. Guadalajara, México.}
\email{alonso.castillor@academicos.udg.mx}

\date{October 2026}

\begin{document}
\raggedbottom

\begin{abstract}
Let $G$ be a group and let $A$ be a set with at least two elements. A lazy cellular automaton over $A^G$ acts as the identity except at the cells where a fixed pattern occurs, where it replaces the central symbol by a fixed different symbol. Let $M(G;A)$ be the submonoid of $\mathrm{CA}(G;A)$ generated by the lazy and the invertible cellular automata. We study when $M(G;A) = \mathrm{CA}(G;A)$. Since no lazy cellular automaton is surjective, every surjective element of $M(G;A)$ is invertible; hence $M(G;A) \neq \mathrm{CA}(G;A)$ whenever $A$ is infinite or $G$ has an element of infinite order. More unexpectedly, if $|A| = 2$ and $G$ has a quotient isomorphic to the Klein four-group $C_2 \times C_2$, then $M(G;A) \neq \mathrm{CA}(G;A)$. For locally finite groups, these are the only obstructions: we prove that $M(G;A) = \mathrm{CA}(G;A)$ if and only if $3 \leq |A| < \infty$, or $|A| = 2$ and $G$ has no quotient isomorphic to $C_2 \times C_2$.
\end{abstract}

\maketitle
\enlargethispage{4pt}

\section{Introduction}

Let $G$ be a group and let $A$ be a set with at least two elements. Cellular automata (CA) are transformations of the configuration space $A^G$ that can be defined via a local map that only depends on a finite neighborhood. The set of all cellular automata over $A^G$ forms a monoid $\CA(G;A)$ under composition, whose group of units is denoted by $\ICA(G;A)$. The structure and the generation of this monoid have been studied in \cite{CeccheriniCoornaert2023,CastilloGadouleau2019,CastilloMinimalGenerators2022,CastilloGenerating2022}.

Lazy cellular automata were introduced in \cite{CastilloIdempotent2024} and further studied in \cite{CastilloUnique2025, AlcalaOrder2026}. They are the CA whose local defining map has a unique active transition $p \in A^S$, where $S$ is a finite subset of $G$ that contains the group identity $e$. Thus, a lazy CA acts as the identity on a configuration except at the coordinates where the pattern $p$ appears; at each of these coordinates, it replaces the central symbol $p(e)$ by a fixed symbol $a\neq p(e)$.

When $S = \{e\}$, lazy CA with unique active transition in $A^{\{e\}}$ may be identified with idempotents of defect $1$ in the full transformation monoid $\text{Tran}(A)$. When $A$ is finite, it is well known \cite{Howie1995} that $\text{Tran}(A)$ is generated by its invertible elements together with its idempotents of defect $1$. This fails for infinite $A$, since a product involving an idempotent of defect $1$ is never surjective; Corollary~\ref{cor:infinite-obstruction} below is the analogue of this fact for cellular automata. This inspired the following problem proposed in \cite[Problem 2]{AlcalaOrder2026}.

\begin{problem}\label{prob:generation}
If $\mathcal{L}(G;A)$ is the set of all lazy cellular automata over $A^G$, prove or disprove the following:
\[ \CA(G;A) = \langle \mathcal{L}(G;A) \cup \ICA(G;A) \rangle.  \]
\end{problem}

In this paper, we answer Problem~\ref{prob:generation} for several classes of groups. Put
\[   M(G;A):=\langle\cL(G;A)\cup\ICA(G;A)\rangle \]
for the submonoid generated by the lazy and invertible cellular automata.
Our main result is the following.

\begin{theorem}\label{thm:main}
Let $G$ be a group and let $A$ be a set with $|A| \geq 2$.
\begin{enumerate}
\item If $A$ is infinite, or if $G$ has an element of infinite order, then $M(G;A) \neq \CA(G;A)$.
\item If $|A| = 2$ and $G$ has a quotient isomorphic to $C_2 \times C_2$, then $M(G;A) \neq \CA(G;A)$.
\item If $G$ is locally finite, then $M(G;A) = \CA(G;A)$ if and only if $3 \leq |A| < \infty$, or $|A| = 2$ and $G$ has no quotient isomorphic to $C_2 \times C_2$.
\end{enumerate}
\end{theorem}

Parts~(1) and~(2) are Corollaries~\ref{cor:infinite-obstruction} and~\ref{cor:klein-quotient}, and part~(3) is Theorem~\ref{thm:locally-finite-classification}. Its special case for finite groups is Theorem~\ref{thm:finite-generation}, which is used in the proof of the general case.

The proof of part~(1) uses the fact that no lazy cellular automaton is surjective; hence every surjective element of $M(G;A)$ is invertible, and it suffices to exhibit a surjective cellular automaton that is not injective. When $G$ and $A$ are both finite, we use the notion of elementary collapsings that map one orbit onto another studied in \cite{CastilloRuiz2023} to obtain generating sets for $\CA(G;A)$, which in this case is equal to $\End_G(A^G)$, the monoid of all $G$-equivariant self-maps of $A^G$. If $|A| = 2$, an unexpected obstruction appears exactly when $G$ has a quotient isomorphic to the Klein-four group $C_2 \times C_2$. Using the tools of restriction and induction of CA developed in \cite{CeccheriniCoornaertInduction2009}, we extend the previous result to locally finite groups. The remaining open cases are those in which $G$ is a periodic group that is not locally finite and either $3 \leq |A| < \infty$, or $|A| = 2$ and $G$ has no Klein-four quotient. 

The structure of the paper is as follows. In Section \ref{sec:preliminaries}, we review the basic notions of lazy CA, and we prove that no lazy CA is injective or surjective, for any group $G$ and any alphabet $A$. In Section \ref{sec:generation}, we prove parts (1), (2), and (3) of Theorem \ref{thm:main}. Finally, in Section~\ref{sec:questions} we discuss some open questions.  


\section{Preliminaries}\label{sec:preliminaries}

In this section, we review standard results and notation. The basic theory of cellular automata over groups may be found in \cite{CeccheriniCoornaert2023}, while the basic theory of lazy cellular automata may be found in \cite{AlcalaOrder2026, CastilloUnique2025, CastilloIdempotent2024}.

Let $G$ be a group with identity element $e$, and let $A$ be a set with
at least two elements. The \emph{full shift}, or \emph{configuration space}, $A^G$ is the set of all functions $x\colon G\to A$. A \emph{pattern} $p \in A^S$ is a function $p : S \to A$, where $S$ is a finite subset of $G$. We use the following left action of $G$ on
$A^G$:
\[
  (g\cdot x)(h)=x(hg),\qquad g,h\in G,\quad x\in A^G.
\]
Given $g\in G, x\in A^G$ and a pattern $p\in A^S$, we say that $p$ \emph{occurs}, or \emph{appears}, at $g$ in $x$ if
\[
    (g\cdot x)\vert_S = p.
\]

A map $\tau\colon A^G\to A^G$ is a cellular automaton if there exist a finite subset $S\subseteq G$, called a \emph{memory set}, or a \emph{neighborhood}, for $\tau$, and a local map $\mu\colon A^S\to A$ such that
\[  \tau(x)(g)=\mu\bigl((g\cdot x)|_S\bigr),  \qquad \forall \ x\in A^G, g\in G. \]
The set of all such maps, equipped with composition, forms a monoid denoted by $\CA(G;A)$. We denote its group of units by $\ICA(G;A)$. We compose maps from right to left, so that $\tau\sigma=\tau\circ\sigma$. An elementary property of cellular automata is that they are always \emph{$G$-equivariant} in the sense that 
\[ \tau(g \cdot x) = g \cdot \tau(x), \quad \forall g \in G, x \in A^G.  \]

\begin{definition}\label{def-lazy}
A cellular automaton $\tau : A^G \to A^G$ is called \emph{lazy} if there exist a finite set $S\subseteq G$ with $e \in S$, a pattern $p \in A^S$, and a symbol $a \in A \setminus \{p(e)\}$ such that
\[    \tau(x)(g)=
  \begin{cases}
    a,    & (g\cdot x)|_S=p,\\
    x(g), & (g\cdot x)|_S \neq p
  \end{cases} \]
  for all $x \in A^G$, $g \in G$. We call $S$ the \emph{active neighborhood} of $\tau$, $p \in A^S$ the \emph{unique active transition} of $\tau$, and $a \in A \setminus \{p(e)\}$ the \emph{writing symbol} of $\tau$. We also denote this by $\tau = \tau_p^a$. Write
$\cL(G;A)$ for the set of lazy cellular automata over $A^G$.
\end{definition}

For $c\in A$, denote
\[  S_c:=p^{-1}(\{c\})\subseteq S,  \qquad   \supp_c(x):=\{g\in G:x(g)=c\}. \]

We shall use the following results, which were proved in
\cite[Lemmas~2 and~3]{AlcalaOrder2026}. Lemma~\ref{lem:new-occurrence} also
appears in the study of idempotent cellular automata
\cite[Lemma~3]{CastilloIdempotent2024} and, in the one-dimensional setting,
in
\cite[Lemma~3.2]{CastilloUnique2025}.

\begin{lemma}\label{lem:new-occurrence}
Let $\tau : A^G \to A^G$ be a lazy cellular automaton with unique active transition $p \in A^S$ and writing symbol $a \in A \setminus \{p(e)\}$. If $\tau(x)|_S=p$ for some $x\in A^G$, then there exists $s\in S_a$ such that $(s\cdot x)|_S=p$.
\end{lemma}

Intuitively, Lemma \ref{lem:new-occurrence} tells us that the occurrence of $p$ in the image $\tau(x)$ forces the occurrence of $p$ in $x$ shifted by an element of $S_a$. 

\begin{lemma}\label{lem:supports}
Let $\tau : A^G \to A^G$ be a lazy cellular automaton with unique active transition $p \in A^S$ and writing symbol $a \in A \setminus \{p(e)\}$. Then, for
every $x\in A^G$ and $i, j\in\mathbb{Z}_{\geq 0}, i\leq j$
\[
  \supp_a(\tau^i(x))\subseteq\supp_a(\tau^j(x)).
\]
Furthermore, if $\supp_a(x)=\supp_a(\tau(x))$, then $\tau(x)=x$.
\end{lemma}

The following elementary results are new.

\begin{proposition}\label{thm:not-injective}
No lazy cellular automaton is injective.
\end{proposition}
\begin{proof}
Let $\tau : A^G \to A^G$ be a lazy cellular automaton with unique active transition $p \in A^S$ and writing symbol $a \in A \setminus \{p(e)\}$. Define $x_0\in A^G$ by
\[
  x_0(g)=
  \begin{cases}
    p(g), & g\in S,\\
    a,    & g\notin S.
  \end{cases}
\]
Put $x_{n+1}=\tau(x_n)$ for $n\geq0$, and let
\[
  C(x):=G\setminus\supp_a(x).
\]
By Lemma~\ref{lem:supports},
\[
  C(x_0)\supseteq C(x_1)\supseteq C(x_2)\supseteq\cdots.
\]
The set $C(x_0)$ is contained in $S$, which is finite, so this chain stabilizes. Let $k$ be the least nonnegative integer such that $C(x_k)=C(x_{k+1})$. Then, Lemma~\ref{lem:supports} shows that
$x_k=x_{k+1}$. Since $x_0|_S=p$, we have $x_1(e) = \tau(x_0)(e) = a \neq p(e)$. Hence,
$x_0\neq x_1$ and $k\geq1$. By the minimality of $k$,
$x_{k-1}\neq x_k$, whereas
\[
  \tau(x_{k-1}) = x_k = x_{k+1} = \tau(x_k).
\]
Therefore, $\tau$ is not injective.
\end{proof}

\begin{proposition}\label{thm:not-surjective}
No lazy cellular automaton is surjective.
\end{proposition}

\begin{proof}
Let $\tau : A^G \to A^G$ be a lazy cellular automaton with unique active transition $p \in A^S$ and writing symbol $a \in A \setminus \{p(e)\}$. Let
\[  F:=S\cup S_aS_a,  \textup{ where } S_aS_a:=\{st:s,t\in S_a\}. \]
Define $q \in  A^F$ by $q|_S =p$ and $q(g) = p(e)$ for $g\in F\setminus S$. We will show that the pattern $q$ does not occur in the image set $\tau(A^G)$, which implies that $\tau(A^G) \neq A^G$.  

Suppose that $\tau(x)|_F=q$ for some $x\in A^G$. Since $\tau(x)|_S=p$, Lemma~\ref{lem:new-occurrence} gives $t\in S_a$ such that $(t\cdot x)|_S=p$. Then, for every $s\in S_a$, we have
\[ a = p(s) = (t\cdot x)(s) = x(st). \] 
Thus, $\tau(x)(st)=a$ for every $s\in S_a$ by Definition \ref{def-lazy}. Since $st\in S_a S_a \subseteq F$, it follows by the hypothesis that $q(st) = \tau(x)(st)  =a$, for every $s\in S_a$. The definition of $q$
gives $q^{-1}(\{a\})=S_a$, and hence $S_at\subseteq S_a$.

Since right multiplication by $t$ preserves cardinality and $S_a$ is a finite subset of $G$, we have $S_at=S_a$. As a consequence of $t\in S_a=S_at$, we have $t=st$ for some $s\in S_a$, so $s=e$. This contradicts $p(e)\neq a$, and the result follows.
\end{proof}


\section{Generation by lazy and invertible cellular automata}
\label{sec:generation}

Recall that
\[
  M(G;A):=\bigl\langle \cL(G;A)\cup\ICA(G;A)\bigr\rangle
  \leq\CA(G;A).
\]
We first give a general obstruction to the equality
$M(G;A)=\CA(G;A)$.

\begin{proposition}\label{prop:surjective-obstruction}
Every surjective element of $M(G;A)$ is invertible.
\end{proposition}

\begin{proof}
Let $\tau=\tau_1\tau_2\cdots\tau_n\in M(G;A)$ be surjective, where every
$\tau_i$ is lazy or invertible. Suppose that at least one factor is lazy,
and let $i$ be the smallest index such that $\tau_i$ is lazy. Then
$\tau_1\cdots\tau_{i-1}$ is bijective, while $\im(\tau_i)$ is a
proper subset of $A^G$ by
Proposition \ref{thm:not-surjective}. Hence,
\[
   \im(\tau) = \im(\tau_1\cdots\tau_n)
  \subseteq
  (\tau_1\cdots\tau_{i-1})\bigl(\im(\tau_i)\bigr)
  \subsetneq A^G.
\]
This contradicts the surjectivity of $\tau$, so $\tau$ has no lazy factor. Therefore, all factors of $\tau$ are invertible, so $\tau$ itself is invertible. 
\end{proof}

\begin{corollary}\label{cor:infinite-obstruction}
Suppose that either $A$ is infinite or $G$ contains an element of infinite
order.  Then
\[
  M(G;A)\neq\CA(G;A).
\]
\end{corollary}

\begin{proof}
By Proposition~\ref{prop:surjective-obstruction}, it is enough to find a
surjective non-injective cellular automaton.

If $A$ is infinite, choose a surjective non-injective map $f\colon A\to A$.
The coordinatewise cellular automaton $f^G : A^G \to A^G$ given by 
\[  f^G(x)(g) :=f(x(g)), \quad \forall g \in G, x \in A^G,  \]
is surjective and non-injective. If $A$ is finite and $G$ has an element of infinite order, then \cite[Exercise~3.24]{CeccheriniCoornaert2023} provides a surjective non-injective cellular automaton $\tau\colon A^G\to A^G$; see \cite[Ch.~3]{CeccheriniCoornaertExercises2023} for a solution.
\end{proof}

\begin{proposition}\label{prop:quotient-descent}
For any normal subgroup $N \trianglelefteq G$ there is a surjective monoid homomorphism
\[  \Phi_N\colon\CA(G;A)\longrightarrow\CA(G/N;A) \]
that sends invertible cellular automata to invertible cellular automata and lazy cellular automata to lazy cellular automata or to the identity. In particular,
\[  M(G;A)=\CA(G;A)\quad\Longrightarrow\quad M(G/N;A)=\CA(G/N;A).\]
\end{proposition}

\begin{proof}
The existence of the surjective monoid homomorphism $\Phi_N\colon\CA(G;A)\longrightarrow\CA(G/N;A)$ is given by \cite[Proposition~1.6.2]{CeccheriniCoornaert2023} as
\[ \Phi_N(\tau) := (\pi^*)^{-1} \circ \tau \vert_{\Fix_N(A^G)} \circ \pi^*, \]
where 
\begin{align}
    \Fix_N(A^G):=\{ x \in A^G : n \cdot x = x, \forall \ n \in N\}, \label{eq:fix_of_subgroup}
\end{align}
and $\pi^* : A^{G/N} \to \Fix_N(A^G)$ is a bijection given by 
\[\pi^*(y)=y \circ \pi, \qquad \forall \ y \in A^{G/N},\]
where $\pi : G \to G/N$ is the canonical projection.

Now let $\tau=\tau_p^a\in\cL(G;A)$, with active neighborhood $S$. Assume $p(s)\neq p(t)$ for some $s,t\in S$ with $\pi(s)=\pi(t)$. Then, $t^{-1}s\in N$, which implies that $t^{-1}s\cdot x = x$, for all $x\in \Fix_N(A^G)$. Thus, 
\[x(s) = x(t), \quad x\in\Fix_N(A^G).\]
Since $p(s)\neq p(t)$, we get that $p$ cannot occur at $e$ in any $N$-fixed configuration. It is known that $\Fix_N(A^G)$ is $G$-invariant (see \cite[Proposition~1.3.6]{CeccheriniCoornaert2023}), so $p$ does not occur in any $N$-fixed configuration. Then, by Definition \ref{def-lazy},
\[\Phi_N(\tau) = \id.\]
On the other hand, if for every $s, t\in S$, $p(s)\neq p(t)$ implies $\pi(s)\neq\pi(t)$, we get that $p$ defines a pattern $\overline p\in A^{\pi(S)}$ by $\overline p(\pi(s))=p(s)$. For $x\in A^{G/N}$, $g\in G$, and $s\in S$, notice that
\[ (g\cdot\pi^*(x))(s) = \pi^\ast(x)(sg) = (x\circ \pi)(sg) = x (\pi(s)\pi(g)) = (\pi(g)\cdot x)(\pi(s)).
\]
Hence, $p$ occurs at $g$ in $\pi^*(x)$ if and only if $\overline p$ occurs at $\pi(g)$ in $x$. Thus
\[\tau_p^a\circ\pi^\ast = \pi^\ast\circ\tau_{\overline{p}}^a.\]
Then, after applying $(\pi^\ast)^{-1}$ we get
\[(\pi^\ast)^{-1}\circ\tau_p^a\circ\pi^\ast = \tau_{\overline{p}}^a.\]
As a consequence, 
\[\Phi_N(\tau_p^a) = \tau_{\overline{p}}^a.\]
Notice that $\tau_{\overline{p}}^a\in\mathcal{L}(G/N; A)$, since  $\overline{p}(\pi(e)) = p(e) \neq a$. Consequently, $\Phi_N$ maps lazy CAs to lazy CAs or the identity. Since $\Phi_N$ is a monoid homomorphism, it follows $\Phi_N(\ICA(G;A)) \subseteq \ICA(G/N;A)$. Therefore,
\[\Phi_N(M(G;A))\subseteq M(G/N;A).\] 
If $M(G;A)=\CA(G;A)$, the surjectivity of $\Phi_N$ gives $\CA(G/N;A)=\Phi_N(M(G;A))\subseteq M(G/N;A)$.
\end{proof}

\subsection{Finite groups}

For a finite group $G$ and any alphabet $A$, we have $\CA(G;A)=\operatorname{End}_G(A^G)$, the monoid of $G$-equivariant transformations of $A^G$, because every $G$-equivariant map $\tau\colon A^G\to A^G$ admits $G$ as a memory set, with local map $\mu(x)=\tau(x)(e)$. When $A$ is finite, this identifies the problem with the generation of the endomorphism monoid of a finite $G$-set (see \cite{CastilloRuiz2023}). 

We first give the required definitions for a finite $G$-set $X$. Write $\operatorname{End}_G(X)$ for the monoid of $G$-equivariant transformations of $X$, and $\operatorname{Aut}_G(X)$ for its group of units. For $x\in X$, denote its orbit and
stabilizer by
\[
  Gx:=\{g\cdot x:g\in G\},
  \qquad G_x:=\{g\in G:g\cdot x=x\}.
\]
For $f\in\operatorname{End}_G(X)$, its \emph{kernel} means the equivalence
relation
\[
  \ker(f):=\{(z,z')\in X\times X:f(z)=f(z')\}.
\]
Equivariance gives
$G_z\leq G_{f(z)}$; equality holds if and only if $f$ is injective on $Gz$.

\begin{definition}\label{def:elementary-collapse}
Let $X$ be a $G$-set. Consider $x, y\in X$  such that $Gx\neq Gy$ and $G_x\leq G_y$. Define
$[x\mapsto y]\colon X\to X$ by
\[
  [x\mapsto y](z)=
  \begin{cases}
    g\cdot y,&z=g\cdot x\text{ for some }g\in G,\\
    z,&z\notin Gx.
  \end{cases}
\]
We call $[x\mapsto y]$ an \emph{elementary collapsing}.
\end{definition}

The inclusion $G_x\leq G_y$ makes $[x \mapsto y]$ well-defined.
It is equivariant and idempotent, since it maps $Gx$ onto the distinct
orbit $Gy$ and fixes $X\setminus Gx$. We call $[x\mapsto y]$
\emph{diagonal} if $G_x=G_y$; this corresponds to the diagonal case of
\cite[Definition~15]{CastilloRuiz2023}. In that case $[x\mapsto y]$ is
injective on $Gx$, so
\begin{equation}\label{eq:diagonal-kernel}
  \ker([x\mapsto y])=\Delta_X\cup
  \{(g\cdot x,g\cdot y),(g\cdot y,g\cdot x):g\in G\},
\end{equation}
where $\Delta_X:=\{(z,z):z\in X\}$. 

Similarly to \eqref{eq:fix_of_subgroup}, given a $G$-set X, and $H\leq G$, define
\[\Fix_H(X) := \{x\in X : H\leq G_x\}.\]

We shall use the following classical result of Burnside; see \cite[Theorem~2.3.2, p.744]{Hazewinkel2000}.
\begin{theorem}\label{thm:Burnside}
    Let $G$ be a finite group, and let $X$ and $Y$ be finite $G$-sets. Then, the following are equivalent:
    \begin{enumerate}
        \item There exists a $G$-equivariant bijection from $X$ to $Y$;
        \item $|\Fix_H(X)| = |\Fix_H(Y)|$, for all $H\leq G$.
    \end{enumerate}
\end{theorem}

\begin{lemma}\label{lem:G-equiv.extension}
    Assume $G$ is a finite group and $X$ is a finite $G$-set. Let $U, V\subseteq X$ be $G$-invariant subsets and $\phi: U\to V$ be a $G$-equivariant bijection. Then, $\phi$ extends to an element of $\operatorname{Aut}_G(X)$.
\end{lemma}
\begin{proof}
    By $G$-equivariance of $\phi$, we get that $\phi$ restricts to a bijection $\Fix_H(U)\to \Fix_H(V)$, for all $H\leq G$. Therefore, it follows for every $H\leq G$ that
    \begin{align*}
        \left|\Fix_H(X\backslash U)\right| &= \left|\Fix_H(X)\right| - \left|\Fix_H(U)\right| \\
        &= \left|\Fix_H(X)\right| - \left|\Fix_H(V)\right| \\
        &= \left|\Fix_H(X\backslash V)\right|.
    \end{align*}
    Since $U$ and $V$ are $G$-invariant, so are $X\backslash U$ and $X\backslash V$. Thus, by Theorem \ref{thm:Burnside} we get that there exists a $G$-equivariant bijection
    \[\phi': X\setminus U \to X\setminus V.\]
    Define $\Phi: X \to X$ as follows:
    \[\Phi(x) := \begin{cases}
        \phi(x), & x\in U \\
        \phi'(x), & x\in X\setminus U
    \end{cases}, \qquad \forall \ x\in X.\]
    Since $\phi$ and $\phi'$ are $G$-equivariant bijections, it follows that $\Phi\in \operatorname{Aut}_G(X)$, such that $\Phi\vert_U = \phi$.
\end{proof}

We use the following constructive generation result and its conjugation
consequence.

\begin{lemma}\label{lem:elementary-generation}
Let a finite group $G$ act on a finite set $X$.
\begin{enumerate}[label=\textup{(\roman*)}]
  \item The elementary collapsings, together with
  $\operatorname{Aut}_G(X)$, generate $\operatorname{End}_G(X)$.
  \item Suppose that $Gx\neq Gy$, $Gp\neq Gq$, and
  $G_x=G_p\leq G_y=G_q$. Then $[x\mapsto y]$ and
  $[p\mapsto q]$ are conjugate by an element of
  $\operatorname{Aut}_G(X)$.
\end{enumerate}
\end{lemma}

\begin{proof}
Part \textup{(i)} is \cite[Lemma~10]{CastilloRuiz2023}. For \textup{(ii)}, define $\phi: Gx\sqcup Gy\to Gp\sqcup Gq$ as
\begin{align*}
    \phi(g\cdot x) = g\cdot p, \qquad\forall \  g \in G, \\
    \phi(g\cdot y) = g\cdot q, \qquad\forall \  g \in G.
\end{align*}
Notice that $\phi$ is well-defined since $Gx\neq Gy$, $Gp\neq Gq$, and $G_x=G_p\leq G_y=G_q$. Furthermore, $\phi$ is a $G$-equivariant bijection as a consequence of $G_x = G_p$ and $G_y = G_q$. Observe that $Gx\sqcup Gy$ and $Gp\sqcup Gq$ are $G$-invariant sets. Thus, by Lemma \ref{lem:G-equiv.extension} there exists $\Phi\in\operatorname{Aut}_G(X)$, such that $\Phi\vert_{(Gx\sqcup Gy)} = \phi$. If $z\in X\setminus Gx$, then $\Phi(z)\in X\setminus Gp$. As a consequence,
\begin{align*}
    (\Phi^{-1}\circ[p\mapsto q]\circ\Phi)(z) &= \Phi^{-1}([p\mapsto q](\Phi(z))) = \Phi^{-1}(\Phi(z)) = z.
\end{align*}
Now, notice that
\begin{align*}
    (\Phi^{-1}\circ[p\mapsto q]\circ\Phi)(x) &= \Phi^{-1}([p\mapsto q](\Phi(x))) = \Phi^{-1}([p\mapsto q](p)) = \Phi^{-1}(q) = y.
\end{align*}
The result follows by $G$-equivariance of $\Phi$ and $[p\mapsto q]$.
\end{proof}

The next lemma produces the elementary collapsings that we shall use.

\begin{lemma}\label{lem:lazy-elementary}
Let $G$ be finite, let $p\in A^G$, let $a\in A\setminus\{p(e)\}$, and let
$\tau=\tau_p^a$ be the lazy cellular automaton with active neighborhood $G$.
Then $\tau=[p\mapsto\tau(p)]$.
\end{lemma}
\begin{proof}
Given $g\in G$, notice that $(g\cdot p)\vert_G = p$ if and only if $g\cdot p = p$. Therefore, $p$ occurs at $g$ in $p$ if and only if $g\in G_p$. As a consequence,
\[
  \tau(p)(g)=
  \begin{cases}
    a, & g\in G_p,\\
    p(g), & g\notin G_p.
  \end{cases}
\]
Assume $g\in G$ satisfies $\tau(p)(g) = a$. Then, $g\in G_p$ or $p(g) = a$ and $g\notin G_p$. Thus
\[\tau(p)^{-1}(\{a\}) = G_p\cup(S_a\setminus G_p).\]
For any $g\in G_p$, notice that 
\[p(g)=(g\cdot p)(e)=p(e)\neq a.\] 
Therefore, $S_a\cap G_p = \emptyset$. This implies
\[|\tau(p)^{-1}(\{a\})| = |G_p| + |S_a| = |G_p| + |p^{-1}(\{a\})|.\]
Now, given $h\in G$, observe the following
\begin{align*}
    g\in p^{-1}(\{a\})h^{-1} \iff p(gh) = a \iff (h\cdot p)(g) = a \iff g\in(h\cdot p)^{-1}(\{a\}).
\end{align*}
Therefore, $(h\cdot p)^{-1}(\{a\})=p^{-1}(\{a\})\,h^{-1}$ for all $h\in G$. Since multiplying by an element of the group preserves the cardinality, we get 
\[|p^{-1}(\{a\})h^{-1}| = |p^{-1}(\{a\})|.\]
It follows that
\[|\tau(p)^{-1}(\{a\})|=|G_p| + |p^{-1}(\{a\})|>|p^{-1}(\{a\})|=|(h\cdot p)^{-1}(\{a\})|, \qquad\forall \  h\in G.
\]
Therefore, $\tau(p)\neq h\cdot p$, for all $h\in G$, which implies $\tau(p)\notin Gp$. It follows that
\[Gp\neq G\tau(p).\]
Now, if $g\in G_p$, by $G$-equivariance of $\tau$ we get that
\[g\cdot \tau(p) = \tau(g\cdot p) = \tau(p).\]
Thus, $G_p\leq G_{\tau(p)}$. This implies that $[p\mapsto\tau(p)]$ is an elementary collapsing. By construction, notice that $\tau(p) = [p\mapsto \tau(p)](p)$. So, by $G$-equivariance we get
\[\tau(g\cdot p) = [p\mapsto \tau(p)](g\cdot p), \qquad \forall \ g\in G.\]
If $x\notin Gp$, then $g\cdot x\neq p$ for all $g\in G$, so $\tau(x)(g) = x(g)$ for all $g\in G$. Therefore,
\[\tau(x) = x = [p\mapsto\tau(p)](x).\]
Consequently, $\tau_p^a = [p\mapsto\tau(p)]$.
\end{proof}

\begin{lemma}\label{lem:diagonal-necessity}
Let $G$ be a finite group, let $X$ be a finite $G$-set and $W\subseteq\operatorname{End}_G(X)$. Assume that $x,y\in X$ satisfy $Gx\neq Gy$ and $G_x=G_y$. If
\[[x\mapsto y]\in \langle W\cup \operatorname{Aut}_G(X)\rangle,\]
then there are $w\in W$ and $\phi, \psi\in\operatorname{Aut}_G(X)$, such that
\[
  [x\mapsto y] = \phi\circ w\circ \psi.
\]
\end{lemma}

\begin{proof}
By hypothesis,
\[
  [x\mapsto y]=\tau_n\tau_{n-1}\cdots\tau_1,
  \qquad
  \tau_i\in W\cup\operatorname{Aut}_G(X).
\]
Notice that $x\neq y$ and $[x\mapsto y](x) = [x\mapsto y](y)$, so, $[x\mapsto y]$ is not invertible. As a consequence, at least one factor is not
invertible. Let $k$ be the least index for which $\tau_k$ is not invertible. Then
$\tau_k\in W$, and
\[
  \psi :=\tau_{k-1}\cdots\tau_1
  \in\operatorname{Aut}_G(X),
\]
where $\psi=\id_X$ if $k=1$. Define $u := \tau_n\cdots \tau_{k+1}$. Again, if $k = n$, we get $u = \id_X$. Notice that
\[
  [x\mapsto y] = u\circ\tau_k\circ \psi.
\]
Consider $v = \tau_k\circ\psi$. On the finite set $X$, an injective self-map is bijective. Also, a bijective equivariant map has an equivariant inverse. Therefore, $\tau_k$ is not injective since it is not invertible. Since $\psi$ is invertible, $v$ is not injective. Now, notice that $[x\mapsto y] = u\circ v$. Then, by \eqref{eq:diagonal-kernel}
\[\ker(v)\subseteq \ker([x\mapsto y]) = \Delta_X\cup
  \{(g\cdot x,g\cdot y),(g\cdot y,g\cdot x):g\in G\}.\]
Since $v$ is not injective, there exist $(z, z')\in \ker(v)$ such that $z\neq z'$. Thus,
\[(z, z')\in \{(g\cdot x,g\cdot y),(g\cdot y,g\cdot x):g\in G\}.\]
Without loss of generality, there exists $g\in G$ such that $g\cdot x = z$ and $g\cdot y = z'$. Then, by $G$-equivariance
\[g\cdot v(x) = v(g\cdot x) = v(g\cdot y) = g\cdot v(y).\]
Thus, for each $h\in G$ it follows that
\[h\cdot v(x) = hg^{-1}g\cdot v(x) =  hg^{-1}g\cdot v(y) = h\cdot v(y).\]
Now, by $G$-equivariance we get
\[v(h\cdot x) = v(h\cdot y), \qquad \forall \ h\in G.\]
As a consequence,
\[\{(g\cdot x,g\cdot y),(g\cdot y,g\cdot x):g\in G\}\subseteq \ker(v).\]
It follows that
\begin{align}\label{eq:ker-equality}
    \ker(v) = \ker([x\mapsto y]).
\end{align}
Define $\varepsilon: v(X)\to [x\mapsto y](X)$, given by
\[\varepsilon(v(z)) = [x\mapsto y](z), \qquad \forall \ z\in X.\]
If $v(z) = v(z')$, then $[x\mapsto y](z) = [x\mapsto y](z')$ by \eqref{eq:ker-equality}. Thus, $\varepsilon$ is well-defined. Notice that for each $[x\mapsto y](z)\in [x\mapsto y](X)$, there exists $v(z)\in v(X)$ such that
\[\varepsilon(v(z)) = [x\mapsto y](z).\]
Therefore, $\varepsilon$ is surjective. Now, assume $\varepsilon(v(z)) = \varepsilon(v(z'))$, for some $z, z'\in X$. Then $[x\mapsto y](z)= [x\mapsto y](z')$. This implies that $(z, z')\in \ker([x\mapsto y])$. So, by \eqref{eq:ker-equality} we get $(z, z')\in\ker(v)$, which implies $v(z)= v(z')$. Thus, $\varepsilon$ is injective. Now, given $z\in X$ and $g\in G$, notice that
\[\varepsilon(g\cdot v(z)) = \varepsilon(v(g\cdot z)) = [x\mapsto y](g\cdot z) = g\cdot[x\mapsto y](z) = g\cdot\varepsilon(v(z)).\]
Therefore, $\varepsilon$ is $G$-equivariant. Notice that $v(X)$ and $[x\mapsto y](X)$ are $G$-invariant by $G$-equivariance of $v$ and $[x\mapsto y]$, respectively. Then, by Lemma \ref{lem:G-equiv.extension} there exists $\phi\in \operatorname{Aut}_G(X)$ that extends $\varepsilon$. For each $z\in X$, notice that
\[\phi(v(z)) = \varepsilon(v(z)) = [x\mapsto y](z).\]
Therefore, $\phi\circ v = [x\mapsto y]$. Then, by construction of $v$ we get
\[[x\mapsto y] = \phi\circ\tau_k\circ\psi,\]
where $\phi, \psi\in \operatorname{Aut}_G(X)$ and $\tau_k\in W$.
\end{proof}

\Needspace{10\baselineskip}
We now consider an obstruction that occurs for a binary alphabet over the Klein group. 

\begin{lemma}\label{lem:lazyCAs-on-binary-Klein}
    Consider $V=C_2\times C_2$ and $A=\{0,1\}$. Define
    \[\mathcal{O}_c := \{x\in A^V : |\supp_{c}(x)| = 1\}, \qquad c\in A.\]
   Notice that $\mathcal{O}_0$ and $\mathcal{O}_1$ are the two $V$-orbits whose stabilizer is $\{e\}$. Then, for every $\tau\in\mathcal{L}(V; A)$ 
\[\tau(\mathcal{O}_0)\cap\tau(\mathcal{O}_1)\neq\emptyset \Longrightarrow \tau(\boldsymbol{0}) = \tau(\boldsymbol{1}), \]
where $\boldsymbol{0}$ and $\boldsymbol{1}$ are the constant configurations in $A^V$ that map $V$ to $0$ and $1$, respectively. 
\end{lemma}
\begin{proof}
    Let $S\subseteq V$ such that $e\in S$. Now, consider $p\in A^S$ and $\tau\in\mathcal{L}(V; A)$ such that $p$ is the unique active transition of $\tau$. We shall prove that $\tau(\boldsymbol{0}) \neq \tau(\boldsymbol{1})$ implies $\tau(\mathcal{O}_0)\cap\tau(\mathcal{O}_1) = \emptyset$. Assume 
    \[\tau(\boldsymbol{0}) \neq \tau(\boldsymbol{1}).\]
    Notice that the writing symbol of $\tau$ is $a\in \{0, 1\}\setminus\{p(e)\}$. Assume $a = 1$. Then, $p(e) = 0$, and 
    \[\tau(x)(g) = x(g), \qquad \forall \ g\in V \textup{ such that } x(g) = 1.\]
    Therefore, $\tau(\boldsymbol{1}) = \boldsymbol{1}$. If $p = \boldsymbol{0}\vert_S$, then $p$ occurs at $g$ in $\boldsymbol{0}$ for each $g\in V$, which implies $\tau(\boldsymbol{0}) = \boldsymbol{1}$, so
    \[\tau(\boldsymbol{0}) = \tau(\boldsymbol{1}), \textup{ contradicting } \ \tau(\boldsymbol{0})\neq \tau(\boldsymbol{1}).\]
    It follows that there is $s_0\in S$ such that
    \[p(s_0) = 1.\]
    Consider $x\in \mathcal{O}_1$ and $h$ the unique element in $V$ such that $x(h) = 1$. If $p$ occurs at $g$ in $x$, for some $g\in V$, then 
    \[(g\cdot x)(s_0) = p(s_0) = 1 \Longrightarrow x(s_0g) = 1. \]
    By uniqueness of $h$, we get that $g = s_0^{-1}h$. Therefore, $p$ occurs at most once in $x$. Then, by Definition \ref{def-lazy} it follows that
    \[\supp_1{\tau(x)}\subseteq\{h, s_0^{-1}h\} \Longrightarrow |\supp_1{\tau(x)}| \leq 2.\]
    If $y\in\mathcal{O}_0$, by Lemma \ref{lem:supports} we get 
    \[\supp_1{y}\subseteq \supp_1{\tau(y)} \Longrightarrow |\supp_1{\tau(y)}|\geq 3.\]
    As a consequence, $\tau(x)\neq \tau(y)$, for each $x\in\mathcal{O}_1, y\in\mathcal{O}_0$, which implies
    \[\tau(\mathcal{O}_0)\cap\tau(\mathcal{O}_1) = \emptyset.\]
    Now, assume that the writing symbol of $\tau$ is $a = 0$. Thus, $p(e) = 1$. Define $\sigma: A^V\to A^V$ as
    \[\sigma(x)(g) = \begin{cases}
        0, & x(g) = 1 \\
        1, & x(g) = 0
    \end{cases} \qquad \forall \  x\in A^V, g\in V.\]
    Notice that $\sigma$ is $V$-equivariant and $\sigma^2 = \id_{A^V}$. Thus, $\sigma\in\operatorname{Aut}_V(A^V)$. Consider $\tau'\in\CA(V; A)$ as
    \[\tau' := \sigma\circ\tau\circ \sigma.\]
    It follows that $\tau'$ is lazy with unique active transition $q\in A^S$ given by
    \[q(s) = \begin{cases}
        0, & p(s) = 1 \\
        1, & p(s) = 0
    \end{cases} \qquad \forall \ s\in S,\]
    and writing symbol $1$. Notice that $\sigma(\mathcal{O}_b) = \mathcal{O}_c$ and $\sigma(\boldsymbol{b}) = \boldsymbol{c}$, for each $b\in A$ and $c\in A\setminus\{b\}$. Therefore, $\tau(\boldsymbol{0})\neq\tau(\boldsymbol{1})$ implies $\tau'(\boldsymbol{0})\neq\tau'(\boldsymbol{1})$. Analogously to the case when the writing symbol of $\tau$ is 1, we get that 
    \[\tau'(\boldsymbol{0}) \neq \tau'(\boldsymbol{1})\Longrightarrow \tau'(\mathcal{O}_0)\cap\tau'(\mathcal{O}_1) = \emptyset.\]
    Thus, by conjugating by $\sigma$ it follows that $\tau(\boldsymbol{0}) \neq \tau(\boldsymbol{1})$ implies $\tau(\mathcal{O}_0)\cap\tau(\mathcal{O}_1) = \emptyset$.
\end{proof}
\begin{corollary}\label{cor:klein-obstruction}
    Let $V=C_2\times C_2$ and $A=\{0,1\}$. Then, $M(V; A)\neq \CA(V; A)$. Furthermore, $[x\mapsto y] \notin M(V; A)$, with $x\in\mathcal{O}_1$ and $y\in\mathcal{O}_0$.
\end{corollary}
\begin{proof}
Let $\varepsilon := [ x \mapsto y]$. Since $A$ and $V$ are both finite, $\varepsilon\in \CA(V; A)$. Now, suppose $\varepsilon\in M(V; A)$. Then, by Lemma \ref{lem:diagonal-necessity} there are $\phi, \psi\in \operatorname{ICA}(V; A)$ and $\tau\in \mathcal{L}(V; A)$, such that
    \[\varepsilon = \phi\circ\tau\circ\psi.\]
    By construction, there are $x\in\mathcal{O}_1$ and $y\in\mathcal{O}_0$ such that $\varepsilon(x) = y$ and $\varepsilon(y) = y$. Then
    \[\phi(\tau(\psi(x))) = \phi(\tau(\psi(y))).\]
    By injectivity of $\phi$ we get $\tau(\psi(x)) = \tau(\psi(y))$. It is known that invertible equivariant maps preserve stabilizers. Since $A^V$ only has two $V$-orbits whose stabilizers are $\{e\}$, it follows that $\psi$ permutes $\mathcal{O}_1$ and $\mathcal{O}_0$, or, it fixes them. In any case, $\tau(\psi(x)) = \tau(\psi(y))$ implies 
    \[\tau(\mathcal{O}_0)\cap\tau(\mathcal{O}_1) \neq \emptyset.\]
    Then, by Lemma \ref{lem:lazyCAs-on-binary-Klein} we get that $\tau(\boldsymbol{0}) = \tau(\boldsymbol{1})$. Since the orbits with the constant configurations are the only $V$-orbits whose stabilizer is $V$, it follows that each invertible $V$-equivariant map permutes them or fixes them. Thus, 
    \[\phi(\tau(\psi(\boldsymbol{0}))) = \phi(\tau(\psi(\boldsymbol{1}))).\]
    As a consequence, $\varepsilon(\boldsymbol{0}) = \varepsilon(\boldsymbol{1})$, which contradicts the fact that $\varepsilon$ fixes all $V$-orbits other than $\mathcal{O}_1$. Therefore, $\varepsilon\notin M(V; A)$. This shows that $M(V; A)\neq\CA(V; A)$. 
\end{proof}

\begin{corollary}\label{cor:klein-quotient}
Let $G$ be any group and let $|A|=2$. If $G$ has a quotient
isomorphic to $C_2\times C_2$, then $M(G;A)\neq\CA(G;A)$.\footnote{The case $G=C_2\times C_2$, namely that the elementary collapsing $[\one_{\{e\}}\mapsto\one_{V\setminus\{e\}}]$ does not belong to $M(V;\{0,1\})$, was found independently by Wenlin Zhang and Haobo Ma and communicated to us privately. See also \cite{KleinFour}.}
\end{corollary}

\begin{proof}
Assume $M(G;A) =\CA(G;A)$. Then, if $N\unlhd G$ satisfies $G/N\cong C_2\times C_2$, by Proposition \ref{prop:quotient-descent} we get $M(G/N; A) = \CA(G/N; A)$, which contradicts Corollary \ref{cor:klein-obstruction}. Therefore, $M(G;A)\neq\CA(G;A)$.
\end{proof}

For a subset $B\subseteq G$, denote by $\one_B\in A^G$ its characteristic
configuration, where we fix two symbols $0,1\in A$. By the definition of the shift action it follows that
\[
  G_{\one_B}=\{g\in G:Bg=B\}.
\]
We proceed to characterize when lazy and invertible CAs generate all CAs for finite groups.

\begin{theorem}\label{thm:finite-generation}
Let $G$ be a finite group and $A$ an alphabet with
$\lvert A\rvert\geq2$.  Then
\[
  M(G;A)=\CA(G;A)
\]
if and only if
\[  3\leq\lvert A\rvert<\infty,
  \qquad\text{or}\qquad
  \lvert A\rvert=2
  \ \text{and $G$ has no quotient isomorphic to $C_2\times C_2$}.\]

\end{theorem}

\begin{proof}
If $A$ is infinite, then $M(G;A)\neq\CA(G;A)$ by Corollary~\ref{cor:infinite-obstruction}. If $|A|=2$ and $G$ has a quotient isomorphic to $C_2\times C_2$, then $M(G;A)\neq\CA(G;A)$ by Corollary~\ref{cor:klein-quotient}. This proves necessity.

Conversely, assume that $3\leq |A| < \infty$, or $|A| = 2$ and $G$ has no quotient isomorphic to $C_2\times C_2$. For $p\in A^G$ and $a\in A\setminus\{p(e)\}$, let $\tau_p^a$ have active neighborhood $G$. By Lemma~\ref{lem:lazy-elementary},
\begin{align*}
    \tau_p^a &= [p\mapsto \tau(p)], \\
  \tau(p)(G_p) &= \{a\}, \\
  \tau(p)|_{G\setminus G_p} &= p|_{G\setminus G_p}.
\end{align*}
By Lemma~\ref{lem:elementary-generation}, it suffices to find, for all $H\leq K\leq G$ such that $G_x=H$ and $G_y=K$ for some $x,y\in A^G$ in distinct orbits, such a $\tau := \tau_p^a$ with $G_p=H$ and $G_{\tau(p)}=K$. Fix $0,1\in A$. We proceed by cases:
\begin{itemize}
    \item[-] Assume $H < K$. Let $p=\one_{K\setminus H}$ and writing symbol $a=1$. Clearly $H\leq G_p$. If $g\in G_p$, we get 
    \[(K\setminus H)g = K\setminus H.\]
    Given $k\in K\setminus H$, it follows that $kg\in K\setminus H$. Then, $kg\in K$, which implies $g\in K$. Then
    \[K\setminus (Hg) = (K\setminus H)g = K\setminus H.\]
    Therefore, $g\in H$. As a consequence, $G_p=H$, and $\tau(p)=\one_K$, so $G_{\tau(p)}=K$.
    
    \item[-] Assume $H = K$ and $|A|\geq 3$. Choose pairwise distinct $a,b,c\in A$, and define $p\in A^G$ as follows:
    \[p(g) := \begin{cases}
        b, & g\in H \\
        c, & g\notin H
    \end{cases}, \qquad\forall \ g\in G.\]
    Notice $G_p=H$. Now, consider $a$ as the writing symbol. Then,
     \[\tau(p)(g) = \begin{cases}
        a, & g\in H \\
        c, & g\notin H
    \end{cases}, \qquad\forall \ g\in G.\]
    Thus, $G_{\tau(p)} = H = K$. 
    
    \item[-] Now, assume $H=K$, $A=\{0,1\}$, and $G$ has no quotient isomorphic to $C_2\times C_2$. Once again we proceed by cases:
    \begin{itemize}
        \item[$\cdot$] If $H=G$, let $p=\boldsymbol{0}$ and $a=1$. Then $\tau(p)=\one_G$.
        
        \item[$\cdot$] If $H\ntrianglelefteq G$, choose $g\notin N_G(H)$. Consider $p=\one_{H\cup gH}$ and the writing symbol $a=0$. Notice that $H\leq G_p$. Also, since $e\in H\cup gH$ we get that
        \[G_p\subseteq H\cup gH.\]
        Assume $G_p\neq H$. Then, $|H|<|G_p|\leq2|H|$ and $|H|$ divides $|G_p|$, since $H\leq G_p$. Hence $[G_p:H]=2$ and $G_p=H\cup gH$. Since the index of $H$ in $G_p$ is 2, we get that $H\trianglelefteq G_p$. As a consequence, 
        \[G_p\leq N_G(H).\]
        Notice that $g\in gH\subseteq G_p$, so $g\in N_G(H)$, which contradicts the choice of $g$. Thus $G_p=H$, and $\tau(p)=\one_{gH}$. Finally, notice that for each $k\in G$
        \[(gH)k = gH \iff Hk = H \iff k\in H,\]
        which implies $G_{\tau(p)} = H$.

        \item[$\cdot$] Consider $H < G$ such that $H\lhd G$. If $G/H\cong C_2$, then the configurations with stabilizer $H$ are $\one_H$ and $\one_{G\setminus H}$, which are in the same orbit, so there are not two distinct orbits with stabilizer $H$. Otherwise, if every element of $G/H$ different from the identity has order 2, then for each $g_1H, g_2H\in G/H$ we have
        \[(g_1H)(g_2H) = (g_1g_2H)^{-1} = (g_2^{-1}H)(g_1^{-1}H) = (g_2H)(g_1H).\]
        As a consequence, $G/H$ is abelian. Then, since $G/H$ is finite, it follows that
        \[G/H\cong (C_2)^r, \textup{ for some } r\geq 2.\]
        Thus, there is a surjective group homomorphism $\nu: G/H\to C_2\times C_2$. By composing the surjective group homomorphism $\pi: G\to G/H$ with $\nu$ we get that $G$ has a quotient isomorphic to $C_2\times C_2$, which is a contradiction. Thus, there is $g_0\in G$ such that
        \[(g_0H)^2\neq H.\]
        Define $p = \one_{H\cup g_0H}$ and set $a = 0$ as the writing symbol. Notice that $H\leq G_p$ since right multiplication by elements of $H$ preserves both cosets $H$ and $g_0H$. Now, define $L := \{eH, g_0H\}\subseteq G/H$. Assume that $L(gH) = L$ for some $gH\in G/H$. Since $eH\in L$, it follows that
        \[gH = (eH)(gH)\in L.\]
        Thus, $gH = eH$ or $gH = g_0H$. If $gH = g_0H$, then
        \[L(g_0H) = \{g_0H, (g_0H)^2\} = \{eH, g_0H\},\]
        which implies $(g_0H)^2 = H$, a contradiction. Therefore, the only element of $G/H$ that preserves $L$ is the identity. Thus, the only elements of $G$ that stabilize $p$ are the elements of $H$. As a consequence,
        \[G_p = H.\]
        Then, $\tau(p) = \one_{g_0H}$ and $G_{\tau(p)} = H$.
        
    \end{itemize}
\end{itemize}
In all previous cases, it follows from Lemma~\ref{lem:elementary-generation} that $M(G;A)=\CA(G;A)$.
\end{proof}


\subsection{Locally finite groups}

We shall use the restriction and induction operations for cellular automata to extend Theorem~\ref{thm:finite-generation} to locally finite groups. Let $H\leq G$. For $x\in A^G$ and $g\in G$, define
$x_g\in A^H$ by $x_g(h)=x(hg)$. The \emph{induction} of
$\tau\in\CA(H;A)$ is given by
\[
  \tau^G(x)(g):=\tau(x_g)(e).
\]
It has the same local map and memory set as $\tau$. Therefore, induction preserves laziness. Also, induction preserves composition and invertibility, since by equivariance of
$\tau$ we get that 
\[(\tau^G(x))_g=\tau(x_g),\]
so $(\tau\sigma)^G=\tau^G\sigma^G$, for each $\sigma\in\CA(H; A)$. Conversely, if $\tau\in\CA(G;A)$
has a memory set contained in $H$, its \emph{restriction} $\tau_H$
is defined by
\[
  \tau_H(w):=\tau(\widetilde w)|_H,
\]
where $\widetilde w\in A^G$ is any extension of $w\in A^H$ to $A^G$. These operations give inverse
monoid isomorphisms between $\CA(H;A)$ and the cellular automata over
$G$ admitting a memory set contained in $H$ (see \cite[Proposition~1.7.2]{CeccheriniCoornaert2023} and \cite{CeccheriniCoornaertInduction2009}).

\begin{lemma}\label{lem:exceptional-stabilizers}
Let $G$ be a finite group, let $A=\{0,1\}$, and define
\[
  \mathcal E(G):=\{N\trianglelefteq G:G/N\cong C_2\times C_2\}.
\]
\begin{enumerate}[label=\textup{(\arabic*)}]
  \item Every elementary collapsing $[x\mapsto y]$ of $A^G$ belongs to
        $M(G; A)$, unless $G_x=G_y\in\mathcal E(G)$.
  \item If $f\in\CA(G; A)$ fixes every configuration whose stabilizer
        belongs to $\mathcal E(G)$, then $f\in M(G; A)$.
\end{enumerate}
\end{lemma}

\begin{proof}
For \textup{(1)}, the constructions in the proof of
Theorem~\ref{thm:finite-generation}, together with
Lemma~\ref{lem:elementary-generation}\textup{(ii)}, apply to $G$ and cover
every elementary collapsing except those such that $G_x=G_y\trianglelefteq G$ and $G/G_x\cong C_2^r$. If $r=1$, no such
collapsing exists, and if $r=2$, then $G_x\in\mathcal E(G)$. Suppose that
$r\geq 3$. Consider $Q = G/G_x$ additively and let $\pi\colon G\to Q$ be the
quotient map. Since $Q\cong C_2^r$ and $Q$ is additive, it can be viewed as a vector space over the finite field $\mathbb{F}_2$. Thus, we can choose linearly independent $u,v,w\in Q$. Now, consider the sets
\[
  U=\{0,u,v,w\},
  \qquad
  U^*=\{u,v,w\}.
\]
First note that if $U + t = U$ for some $t \in Q$, then $t \in U$ because $0 \in U$. By linear independence, the only element in $U$ that preserves $U$ under translation is $0$. The same applies to $U^\ast$, since every nontrivial translation is a fixed-point-free involution, so it cannot preserve the set $U^*$ of odd size. Thus the configurations
\[
  p=\one_{\pi^{-1}(U)},
  \qquad
  q=\one_{\pi^{-1}(U^*)}
\]
satisfy $G_p=G_q=G_x$. By Lemma~\ref{lem:lazy-elementary}, $\tau=[p\mapsto q]$, where $\tau$ is the lazy CA whose active neighborhood is $G$, its unique active transition is $p$ and its writing symbol is $0$. Lemma~\ref{lem:elementary-generation}\textup{(ii)}
completes the proof of \textup{(1)}.

For \textup{(2)}, partition $A^G$ into the sets
$\mathcal C_1,\ldots,\mathcal C_m$ of configurations whose stabilizers lie
in a fixed conjugacy class of subgroups, ordered so that the stabilizer orders
are nonincreasing. Assume $f\in \CA(G; A)$ and that $f$ fixes every configuration whose stabilizer belongs to $\mathcal E(G)$. By $G$-equivariance of $f$, $G_z\leq G_{f(z)}$ for every $z\in A^G$. Thus
\[
  f(\mathcal C_i)\subseteq\bigcup_{1\leq j \leq i}\mathcal C_j, \qquad \forall \ i\in\{1, \ldots, m\}
\]
since the stabilizer can increase in size, moving the image to an earlier class. If its size stays the same, inclusion forces equality of stabilizers, so the image remains in $C_i$. Now, for each $i\in \{1, \ldots, m\}$ define $f_i : A^G \to A^G$ as follows 
\[
  f_i(z) :=
  \begin{cases}
    f(z), & z\in\mathcal C_i\\
    z, & z\notin\mathcal C_i
  \end{cases}, \qquad \forall \ z\in A^G.
\]
Notice that each $f_i$ is $G$-equivariant since $C_i$ is $G$-invariant. Now, finiteness of $G$ implies that $f_i\in\CA(G; A)$ for each $i$. Given $z\in\mathcal{C}_i$, notice that
\[f_j(z) = z, \qquad \forall \ 1\leq j \leq i -1.\]
Since $f_i(z) = f(z)$, $f(z)\in \mathcal{C}_{j'}$ for some $j'\leq i$. Thus,
\[f_j(f(z)) = f(z), \qquad \forall \ i+1\leq j \leq m.\]
As a consequence,
\[f = f_m\circ\cdots\circ f_1.\]
If $\mathcal C_i$ corresponds to an element of $\mathcal E(G)$, then
$f_i=\id$ since $f$ fixes the configurations whose stabilizer is in $\mathcal E(G)$. For every other $\mathcal C_i$, the orbit-type factorization
in \cite[proof of Lemma~10]{CastilloRuiz2023} expresses $f_i$ using
automorphisms and elementary collapsings with source orbit contained in
$\mathcal C_i$. All these collapsings belong to $M(G; A)$ by \textup{(1)}. Hence
$f_i\in M(G; A)$ for every $i$. As a consequence, $f\in M(G; A)$.
\end{proof}

\begin{theorem}\label{thm:locally-finite-classification}
Let $G$ be a locally finite group and let $A$ be a set with
$|A|\geq2$. Then
\[
  M(G;A)=\CA(G;A)
\]
if and only if
\[
  3\leq|A|<\infty,
  \qquad\text{or}\qquad
  |A|=2\ \text{and $G$ has no quotient isomorphic to } C_2\times C_2.
\]
\end{theorem}

\begin{proof}
Necessity follows from Corollaries~\ref{cor:infinite-obstruction}
and~\ref{cor:klein-quotient}.

Suppose first that $3\leq|A|<\infty$, and let $\tau\in\CA(G;A)$. Choose
a finite subgroup $H\leq G$ containing a memory set for $\tau$. By
Theorem~\ref{thm:finite-generation}, $\tau_H\in M(H;A)$. Since induction
preserves composition, invertibility, and laziness, we obtain
\[
  \tau=(\tau_H)^G\in M(G;A).
\]

Now assume that $A=\{0,1\}$ and that $G$ has no quotient isomorphic to
$C_2\times C_2$. Write
\[
  G^2:=\langle g^2:g\in G\rangle.
\]
Since $G/G^2$ is an elementary abelian $2$-group,
\[
  |G/G^2|\leq2.
\]

Let $H\leq G$ be finite. Choose $h_0\in H$ such that
\[
  H\subseteq\langle h_0\rangle G^2.
\]
For each $h\in H$, choose $\delta_h\in\{0,1\}$ and
$g_{h,1},\ldots,g_{h,m_h}\in G$ satisfying
\[
  hh_0^{-\delta_h}=g_{h,1}^2\cdots g_{h,m_h}^2.
\]
The subgroup
\[
  K:=\langle H,\ g_{h,j}:h\in H,\ 1\leq j\leq m_h\rangle
\]
is finite, and $H\subseteq\langle h_0\rangle K^2$. Consequently, for every
$N\in\mathcal E(K)$,
\begin{equation}\label{eq:small-quotient-image}
  K^2\leq N,
  \qquad
  [H:H\cap N]=|HN/N|\leq2.
\end{equation}

Let $\tau\in\CA(G;A)$, and choose $H$ to contain a memory set for
$\tau$. By Lemma~\ref{lem:elementary-generation}\textup{(i)}, $\tau_H$ is
a product of automorphisms and elementary collapsings. Automorphisms induce to
invertible cellular automata. By
Lemma~\ref{lem:exceptional-stabilizers}\textup{(1)}, every elementary
collapsing also induces to an element of $M(G;A)$, except possibly one of
the form
\[
  \varepsilon=[p\mapsto q],
  \qquad
  H_p=H_q\in\mathcal E(H).
\]

Choose $K\geq H$ as above. We claim that $\varepsilon^K$ fixes
$\Fix_N(A^K)$ pointwise for every $N\in\mathcal E(K)$. Indeed, for $x\in \Fix_N(A^K), g\in K, n\in H\cap N$, and $h\in H$,
\begin{align*}
    (n\cdot x_g)(h) &= x_g(hn) \\
    &= x(hng) \\
    &= x(hg(g^{-1}ng)) \\
    &= x(hg) \tag{$g^{-1}ng\in N$}\\
    &= x_g(h).
\end{align*}
Thus, $H\cap N\leq H_{x_g}$. Then, by
\eqref{eq:small-quotient-image}
\[|Hx_{g}| = [H: H_{x_g}] \leq [H : H\cap N]\leq 2.\]
Notice that $|Hp|=[H:H_p]=4$ since $H_p\in\mathcal{E}(H)$. Therefore, $x_g\notin Hp$, which implies $\varepsilon(x_g)=x_g$. It follows that
\[
  \varepsilon^K(x)(g)=\varepsilon(x_g)(e)=x(g),
  \qquad \forall \ g\in K,
\]
which proves the claim.

Lemma~\ref{lem:exceptional-stabilizers}\textup{(2)} now gives
$\varepsilon^K\in M(K;A)$. Inducing its factors to $G$, and using
transitivity of induction, yields
\[
  \varepsilon^G=(\varepsilon^K)^G\in M(G;A).
\]
Therefore every factor of $\tau_H$ induces to an element of $M(G;A)$, and
hence $\tau\in M(G;A)$.
\end{proof}

\begin{remark}
The passage from $H$ to $K$ in the proof of Theorem~\ref{thm:locally-finite-classification} cannot in general be avoided. Let $D:=C_{2^\infty}\rtimes\langle s\rangle$, with $s$ acting on the Pr\"ufer $2$-group by inversion. Then $D$ is locally finite, and since $srs^{-1}r^{-1}=r^{-2}$ and squaring is surjective on $C_{2^\infty}$, we have $D/[D,D]\cong C_2$. Hence $D$ has no Klein-four quotient, and $M(D;\{0,1\})=\CA(D;\{0,1\})$. However, every finite subgroup of $D$ containing $s$ and an element of order $4$ is dihedral of order at least $8$, and therefore has a Klein-four quotient.
\end{remark}

\section{Open questions}\label{sec:questions}

We have determined exactly when $M(G;A)=\CA(G;A)$ for every locally finite
group $G$, and extended the binary obstruction to every group with a
Klein-four quotient.

For $|A|=2$, Corollary~\ref{cor:klein-quotient} rules out the equality
$M(G;A)=\CA(G;A)$ for every group $G$ with a Klein-four quotient. This
includes the first Grigorchuk group $\Gamma$, which is finitely generated,
infinite, and periodic, and satisfies $\Gamma/[\Gamma,\Gamma]\cong C_2^3$;
see \cite[Sections~1.6.1 and~7.4]{BartholdiGrigorchukSunic2003}. Thus the
quotient obstruction also applies to periodic groups that are not locally
finite.

There is also an obstruction coming from nonamenability. By
\cite[Theorem~1.2]{Bartholdi2010}, every nonamenable group $G$ admits a
surjective, non-pre-injective cellular automaton over some finite alphabet
$A$. Here non-pre-injectivity means that there exist two distinct
configurations which differ at only finitely many coordinates and have the
same image.
In particular, this cellular automaton is not injective, and
Proposition~\ref{prop:surjective-obstruction} gives
$M(G;A)\neq\CA(G;A)$ for that alphabet. The alphabet is part of the
construction: this conclusion is not asserted for every prescribed finite
alphabet, and in particular not necessarily for a binary one.

\begin{question}
For which pairs $(G,A)$, with $G$ an infinite periodic group that is
not locally finite and $A$ finite, does $M(G;A)=\CA(G;A)$ hold?
\end{question}


\section*{Acknowledgments} 
The first author was supported by SECIHTI \textit{Becas nacionales para estudios de posgrados.}

\section*{Declaration on the use of artificial intelligence}

The authors used ChatGPT 6 Astra and Claude Opus 5.5 to assist with mathematical exploration and discovery. The authors wrote and verified all arguments, and take
full responsibility for the manuscript.


\end{document}